\documentclass[letterpaper, 10 pt, conference]{ieeeconf}  

\IEEEoverridecommandlockouts                              
\usepackage{graphicx}
\usepackage{amssymb, amsmath, bm, cases, comment, color}
\usepackage{subcaption}
\usepackage{algorithm}
\usepackage{algpseudocode}

\newcommand{\tr}{\mathsf{T}}
\newtheorem{prob}{Problem}

\title{\LARGE \bf
Data-driven Koopman Operator-based Prediction and Control \\Using Model Averaging
}

\author{Daisuke Uchida and Karthik Duraisamy
}

\begin{document}

\maketitle
\thispagestyle{empty}
\pagestyle{empty}

\begin{abstract}
This work presents a data-driven Koopman operator-based modeling method using a model averaging technique. 
While the Koopman operator has been used for data-driven modeling and control of nonlinear dynamics, it is challenging to accurately reconstruct unknown dynamics from data and perform different decision-making tasks, mainly due to its infinite dimensionality and difficulty of finding invariant subspaces.
We utilize ideas from a Bayesian inference-based model averaging technique to devise a data-driven method that first populates multiple Koopman models starting 
with a feature extraction using neural networks and then computes point estimates of the posterior of predicted variables.
Although each  model in the ensemble is not likely to be accurate enough for a wide range of operating points or unseen data, the proposed \textit{weighted} linear embedding model combines the outputs of model ensemble aiming at compensating the modeling error of each model so that the overall performance will be improved.
\end{abstract}

\section{INTRODUCTION}
While conventional systems modeling methods describe the dynamics on a state-space, operator theoretic approaches offer an alternative view of the system behaviors through the lens of functions, often called feature maps or observables. 
The Koopman operator, a linear operator acting on a space of feature maps, has been widely used in the context of data-driven modeling, analysis, and control of nonlinear dynamical systems\cite{Koopman_book, Data_driven_book}.
It allows linear evolution of the embedded systems on a function space even if the original dynamics is nonlinear in the state-space. 
Also, it can be approximated numerically from data using either linear regression or neural network training.
The Koopman operator framework is especially appealing when applied to controller synthesis since linear controller designs such as Linear Quadratic Regulator (LQR) and linear Model Predictive Control (MPC) can be utilized with the obtained linear embedding model.

However, obtaining accurate and reliable Koopman operator-based models that can be successfully applied to decision-making tasks is still challenging.
Specifically, it is known that the convergence property of EDMD does not hold for control systems\cite{KORDA_Koopman_MPC}, which implies that Koopman operator-based models may fail to accurately describe control systems even if the training data is sufficiently rich and unbiased. 
To overcome the modeling error of such Koopman operator-based models for control applications, several data-driven controller designs have been proposed that take model uncertainty into account\cite{tube_based_MPC, DeSKO, data-driven_Koopman_H2, handling_plant_model_mismatch}.
Also, if one restricts the attention to control affine systems instead of general nonlinear dynamics, finite-data error bounds are available and stability guarantees may be established on the basis of robust control theories\cite{finite_data_error_bounds, Koopman_control_w_stability_guarantees}.
While many efforts have been made on the model uncertainty of Koopman operator-based models from a controller design perspective, it is also of great importance to devise a modeling method that can yield accurate and generalizable models for control systems.
In \cite{Koopman_control_aware}, a model refinement technique is proposed to incorporate data from closed-loop dynamics into model learning so that the modeling error when applied to controller design problems can be directly mitigated.
To improve the generalizability of the Koopman operator-based models for a wide range of applications, \cite{Koopman_oblique_projection} develops a two-stage learning method utilizing the oblique projection in the context of linear operator learning in a Hilbert space.

In this paper, we adopt a Bayesian approach to learn unknown control systems with the use of the Koopman operator, which takes into consideration that a single Koopman operator-based model is not likely to be accurate and reliable for a wide range of operating points and it will be beneficial to aggregate an ensemble of models and utilize them to find the most useful output.
Ensemble approaches are a popular class of methods that aim to improve the accuracy in learning problems by combining predictions of multiple models\cite{ensemble_learning_survey, ensemble_regression_survey}.
We show that the Bayesian Model Averaging (BMA)\cite{BMA_survey}, a Bayesian inference-based model averaging technique, yields a Koopman operator-based model whose parameters are represented as sums of corresponding parameters of individual models weighted by the posterior model evidence. 
The proposed \textit{Koopman Model Averaging (KMA)} first populates multiple models, starting from a base model whose feature maps are parameterized by a neural network. 
Computing point estimates of the original state and the embedded state then results in the same type of linear embedding model while being model-uncertainty aware.
The outputs of the individual models are incorporated into this uncertainty-aware model
and the overall performance is expected to be improved.

In Section \ref{sec. Koopman operator framework}, the Koopman operator framework for control systems is reviewed. Model learning with neural networks is formulated in Section \ref{sec. learning models from data}, and the proposed modeling approach is derived in Section \ref{sec. weighted Koopman based models}. 
Numerical examples are provided in Section \ref{sec. numerical examples} for evaluation of the proposed method.

\section{Koopman Operator Framework}
\label{sec. Koopman operator framework}
Consider a discrete-time, non-autonomous dynamical system:
\begin{align}
	x_{k+1} = f(x_k, u_k),
	\label{eq. target dynamics}
\end{align}
where $x_k\in \mathcal{X}\subseteq \mathbb{R}^n$, $u_k\in \mathcal{U}\subseteq \mathbb{R}^p$, and $f:\mathcal{X}\times \mathcal{U}\rightarrow \mathcal{X}$ are the state, the input, and the (possibly nonlinear) state-transition mapping, respectively.
It is assumed throughout the paper that the dynamics \eqref{eq. target dynamics} is unknown while we have access to data of the form 
$\{(x_k,u_k,y_k)\mid y_k=f(x_k,u_k)\}$.

We embed the state $x_k$ into a latent space by applying feature maps, or observables, $g:\mathcal{X}\rightarrow \mathbb{R}$.
Let $\mathcal{G}$ denote the function space to which the feature maps $g$ belong.
If the input is constant s.t. $u_k\equiv \bar{u}$, $\forall k$, \eqref{eq. target dynamics} induces autonomous dynamics:
\vspace{-2mm}
\begin{align}
	x_{k+1} = f_{\bar{u}}(x_k):=f(x_k,\bar{u}).
	\label{eq. induced autonomous dynamics}
\end{align}

The state transition through $g$ is then represented by
\begin{align}
	g(x_{k+1})
	=
	(g\circ f_{\bar{u}})(x_k)
	=:(\mathcal{K}_{\bar{u}}g)(x_k),
	\label{eq. Koopman operator for autonomous dynamics}
\end{align}
where the composition operator $\mathcal{K}_{\bar{u}}:\mathcal{G}\rightarrow \mathcal{G}:g\mapsto g\circ f_{\bar{u}}$ is called the Koopman operator associated with the autonomous dynamics \eqref{eq. induced autonomous dynamics}.
It is obvious from \eqref{eq. Koopman operator for autonomous dynamics} that $\mathcal{K}_{\bar{u}}$ is a linear operator and it describes the possibly nonlinear dynamics \eqref{eq. induced autonomous dynamics} linearly in the latent space $\mathcal{G}$.
Since $\mathcal{K}_{\bar{u}}$ is infinite dimensional acting on the function space $\mathcal{G}$, its finite dimensional approximation is often considered for the purpose of dynamical systems modeling.
Given $N_x$ feature maps $g_i\in \mathcal{G}$ ($i=1,\cdots,N_x$), there exists a matrix $K\in \mathbb{R}^{{N_x}\times {N_x}}$ s.t. 
\begin{align}
	[\mathcal{K}_{\bar{u}}g_1\ \cdots \ \mathcal{K}_{\bar{u}}g_{N_x}]^\tr 
	=
	K[g_1\cdots g_{N_x}]^\tr,
	\label{eq. invariance property of finite dimensional approx. autonomous}
\end{align}
if and only if $\text{span}(g_1,\cdots,g_{N_x})$ is an invariant subspace, i.e., $\mathcal{K}_{\bar{u}}g\in \text{span}(g_1,\cdots, g_{N_x})$ for $\forall g\in \text{span}(g_1,\cdots, g_{N_x})$.
Note that finding feature maps $g_i$ that form an invariant subspace may be challenging in practice.
In the data-driven setting, one can approximately obtain $K$ by solving the following linear regression problem:
\begin{align}
	K_\text{approx}:=
	\underset{K}{
		\text{argmin}
	}
	\sum_{i}
	\left\|
	\bm{g}(y_{i}) - K\bm{g}(x_i)
	\right\|_2^2,
\end{align}
where $\bm{g}(x_k):=[g_1(x_k)\cdots g_{N_x}(x_k)]^\tr$ and the training data is given as $\{ (x_i,y_i)\mid y_i=f_{\bar{u}}(x_i) \}$.
It admits the unique analytical solution with the use of pseudo inverse and this procedure is called Extended Dynamic Mode Decomposition (EDMD)\cite{Williams2015}.
Note that the design of the feature maps $g_i$ needs to be pre-specified in EDMD.
The finite dimensional approximation $K_\text{approx}$ then yields a linear embedding model of the form:
\vspace{-1mm}
\begin{align}
	\left\{
	\begin{array}{l}
		g^+ = K_\text{approx} \bm{g}(x_k), 
		\\
		x_{k+1}^\text{pred} = W g^+,
	\end{array}
	\right. 
	\label{eq. linear embedding model autonomous}
\end{align}
where $x_{k+1}^\text{pred}$ is the predicted state at the next time step and the decoder $W\in \mathbb{R}^{n\times {N_x}}$ is similarly obtained by solving linear regression on the given data set.
The decoder can be also introduced as a nonlinear mapping depending on the problem.

For general non-autonomous dynamics \eqref{eq. target dynamics} with control inputs $u_k$, the corresponding Koopman operator can be defined as follows\cite{KORDA_Koopman_MPC}.
For the space of input sequences:
$
	l(\mathcal{U}):=
	\{
		(u_0,u_{1},\cdots)\mid u_k\in \mathcal{U}, \forall k
	\}
$, 
consider a mapping $\hat{f}:\mathcal{X}\times l(\mathcal{U})\rightarrow \mathcal{X}\times l(\mathcal{U}):(x,(u_0,u_1,\cdots))\mapsto (f(x,u_0), (u_1,u_2,\cdots))$.
Also, let $\hat{g}:\mathcal{X}\times l(\mathcal{U})\rightarrow \mathbb{R}$ be an embedding feature map from an extended space $\mathcal{X}\times l(\mathcal{U})$ to $\mathbb{R}$.
Then, the Koopman operator associated with the non-autonomous dynamics \eqref{eq. target dynamics} is defined as a linear operator $\mathcal{K}:\hat{\mathcal{G}}\rightarrow \hat{\mathcal{G}}:\hat{g}\mapsto \hat{g}\circ \hat{f}$ s.t. $\hat{\mathcal{G}}$ is a function space to which feature maps $\hat{g}:\mathcal{X}\times l(\mathcal{U})\rightarrow \mathbb{R}$ belong and the dynamics \eqref{eq. target dynamics} along with a sequence $(u_k,u_{k+1},\cdots)$ of future inputs is represented by
\begin{align}
	\hat{g}(x_{k+1},(u_{k+1},u_{k+2},\cdots))&=
	(\hat{g}\circ \hat{f})(x_k,(u_k,u_{k+1},\cdots))&\nonumber
\\
	&=
	(\mathcal{K}\hat{g})(x_k,(u_k,u_{k+1},\cdots)).
	&\nonumber
\end{align}

The argument on the subspace invariance of the Koopman operator (eq. \eqref{eq. invariance property of finite dimensional approx. autonomous}) also holds for the non-autonomous case, which is written as 
\begin{align}
	[\mathcal{K}\hat{g}_1\cdots \mathcal{K}\hat{g}_{\hat{N}}]^\tr = 
	\hat{K} [\hat{g}_1\cdots \hat{g}_{\hat{N}}],
	\label{eq. invariance property non autonomous} 
\end{align}
where $\hat{K}\in \mathbb{R}^{\hat{N}\times \hat{N}}$ and ${\hat{N}}$ denotes the number of feature maps.
In analogous to \eqref{eq. linear embedding model autonomous}, if we consider $\hat{N}=N_x+p$ feature maps $\hat{g}_i$ ($i=1,\cdots, N_x+p$) of the form:
\begin{align}
	&[\hat{g}_1(x_k,(u_k,u_{k+1},\cdots))\cdots \hat{g}_{N_x+p}(x_k,(u_k,u_{k+1},\cdots))]^\tr
	&\nonumber
\\
	&
	=[g_1(x_k)\cdots g_{N_x}(x_k)\ u_k^\tr]^\tr, 
	&
\end{align}
the first $N_x$ rows of \eqref{eq. invariance property non autonomous} reads
$\bm{g}(x_{k+1}) = \hat{A} \bm{g}(x_k) + \hat{B}u_k$, where $[\hat{A}\ \hat{B}]\in \mathbb{R}^{{N_x}\times (N_x+p)}$ denotes the first $N_x$ rows of $\hat{K}$.
Similar to \eqref{eq. linear embedding model autonomous}, this yields a linear embedding model:
\begin{align}
	\left\{
	\begin{array}{l}
		g^+ = A \bm{g}(x_k) + Bu_k, 
		\\
		x_{k+1}^\text{pred} = C g^+,
	\end{array}
	\right. 
	\label{eq. linear embedding model non autonomous}
\end{align}
in which the parameters $A,B$, and $C$ may be obtained by EDMD with training data $\{ (x_i,u_i,y_i)\mid y_i=f(x_i,u_i) \}$, i.e., 
\begin{align}
	[A\ B]&=
	\underset{[A\ B]}{\text{argmin}}
	\sum_{i}
	\left\|
	\bm{g}(y_{i}) - [A\ B]
	\left[
	\begin{array}{c}
		\bm{g}(x_i)
		\\
		u_i
	\end{array}
	\right]
	\right\|_2^2,&
	\label{eq. EDMD non autonomous 1}
	\\
	C&=
	\underset{C}{\text{argmin}}
	\sum_{i}
	\left\|
	x_i - 
	C\bm{g}(x_{i})
	\right\|_2^2.&
	\label{eq. EDMD non autonomous 2}
\end{align}

A notable feature of the model \eqref{eq. linear embedding model non autonomous} is that the model dynamics in the embedded space is given as a linear time-invariant system and linear controller designs can be utilized to control the possibly nonlinear dynamics \eqref{eq. target dynamics}.
For instance, with the parameters $A$ and $B$ in \eqref{eq. linear embedding model non autonomous}, one can compute an LQR gain that stabilizes the following virtual system:
\begin{align}
	\xi_{k+1} = A\xi_{k+1} + Bu_k,\ \xi_{k}\in \mathbb{R}^{N_x}.
\end{align}

Note that if $\bm{g}$ span an invariant subspace, we have the exact relation: $\bm{g}(x_{k+1})=A\bm{g}(x_k)+Bu_k$.
Since the model \eqref{eq. linear embedding model non autonomous} has a linear decoder, the quadratic loss $\bm{g}(x_k)^\tr Q\bm{g}(x_k)$ of the LQR problem ($Q$ is a weight matrix) can be associated with the original state $x_k$ by $\bm{g}(x_k)^\tr Q\bm{g}(x_k)\approx x_k^\tr(C^\tr QC)x_k$.

While we adopt linear embedding models represented in the form \eqref{eq. linear embedding model non autonomous}, which is a common choice in the literature, more expressive ones such as bilinear models are also available
if the strict linearity w.r.t. $u_k$ in \eqref{eq. linear embedding model non autonomous} is not enough to reconstruct the target dynamics\cite{Koopman_bilinearization, advantages_bilinearization}.  
Also, it is noted that in addition to the extension of the Koopman operator framework to control systems reviewed in this section\cite{KORDA_Koopman_MPC}, there are other formalisms to utilize the Koopman operator for modeling non-autonomous systems\cite{Koopman_control_family}.

\section{Learning Models from Data}
\label{sec. learning models from data}
The accuracy of the linear embedding model \eqref{eq. linear embedding model non autonomous} depends on the design of the feature maps $\bm{g}$ as well as how the model dynamics parameters $A$, $B$, and $C$ are obtained. 
While EDMD provides a simple model learning procedure with the analytic solution, the design of feature maps needs to be user-specified such as monomials and Fourier basis functions and it may not be sufficiently expressive for complex nonlinear dynamics.
Neural networks, on the other hand, are a popular choice of the feature map design since they allow greater expressivity of the model compared to EDMD (e.g., \cite{Physics-based_robabilistic_learning, Learning_Koopman_Invariant_Subspaces,deep_learning_Koopman_CDC2020}). 
With the feature maps $\bm{g}$ characterized by a neural network, both the model dynamics parameters and the feature maps can be learned simultaneously, which is formulated as the following problem:
\begin{prob}
	\vspace{5mm}
	Let $\bm{g}(\cdot;\theta_g):\mathcal{X}\rightarrow \mathbb{R}^{N_x}$ be a neural network characterized by parameters $\theta_g$. 
	Find $\theta_g$, $A\in \mathbb{R}^{N_x\times N_x}$, $B\in \mathbb{R}^{N_x\times p}$, and $C\in \mathbb{R}^{n\times N_x}$ that minimize the loss function:
	\begin{align}
		&J(\theta_g, A,B,C):=
		\sum_{i} 
		\lambda_1
		\left\|
			A\bm{g}(x_i;\theta_g) + Bu_i - \bm{g}(y_i;\theta_g)
		\right\|_2^2
		&\nonumber
	\\
		&
		\hspace{15mm}
		+
		\lambda_2
		\left\|
		C(A\bm{g}(x_i;\theta_g) + Bu_i) - y_i
		\right\|_2^2,&
		\label{eq. loss function}
	\end{align}
	where the data set is given in the form $\{ (x_i,u_i,y_i)\mid y_i=f(x_i,u_i) \}$ and $\lambda_1,\lambda_2\in \mathbb{R}$ are hyperparameters.
	\vspace{5mm}
	\label{problem 1}
\end{prob}

Although neural network-based models are expected to be more expressive and accurate than EDMD-based ones, Problem \ref{problem 1} is typically a high-dimensional non-convex problem and the resulting models may suffer from overfitting or poor learning. 
For instance, solving Problem \ref{problem 1} can lead to inaccurate models if the optimization is terminated at a local minimum with a high loss value, or if the quantity and/or quality of training data are not sufficient to reconstruct the target dynamics for a wide range of operating points.
Therefore, one may need to repeat the data-collection and learning processes multiple times to obtain a satisfying model in practice, which often takes a long time and consumes a large amount of computational resources.

As the first step of the proposed method, we execute Step \ref{problem 1} to obtain a base model.
Considering that the base model may not be perfect, we aim to improve its accuracy further by combining multiple models in the second step. 
Specifically, 
an ensemble of linear embedding models \eqref{eq. linear embedding model non autonomous} is generated using additional data points with the design of feature maps fixed to that of the base model.
These models, including the base model, share the same model structure but are trained on different data sets and their predictive capabilities vary. 
Therefore, even if some model is most accurate for certain unseen data among the model ensemble, others may show better accuracy for different data.
To handle this lack of generalizability of individual models, we use all the models so that the accuracy for unknown regimes of dynamics will be improved.
Specifically, a Bayesian inference-based model averaging technique is utilized to merge the individual models into a new linear embedding model.

\section{Weighted Koopman Operator-based Models}
\label{sec. weighted Koopman based models}
In this section, we first outline the Bayesian Model Averaging (BMA). 
It is then followed by the formulation of the proposed Koopman operator-based model averaging method.
\subsection{Bayesian Model Averaging}
\label{sec. BMA outline}
Given data $\mathcal{D}=\{ (x_i,u_i,y_{i})\mid y_{i}=f(x_i,u_i) \}$ and an ensemble of $N$ models denoted by $\mathbb{M}_i$ ($i=1,\cdots,N$), the posterior distribution of a quantity of interest $q$ is represented as
\vspace{-2mm}
\begin{align}
	p(q\mid \mathcal{D})=
	\sum_{i=1}^{N} p(\mathbb{M}_i\mid \mathcal{D})p(q\mid \mathbb{M}_i,\mathcal{D}),
    \label{eq. BMA posterior}
\end{align}
where
the posterior model evidence is given by
\vspace{-1.5mm}
\begin{align}
	p(\mathbb{M}_i\mid \mathcal{D}) 
	&=
	\cfrac{
		p(\mathcal{D}\mid \mathbb{M}_i)p(\mathbb{M}_i)
	}{
		\sum_{l=1}^{N} p(\mathcal{D}\mid \mathbb{M}_l)p(\mathbb{M}_l)
	}.
	&
\end{align}

The marginal likelihood of model $\mathbb{M}_i$ is represented by
\vspace{-1.5mm}
\begin{align}
	p(\mathcal{D}\mid \mathbb{M}_i)=
	\int p(\mathcal{D}\mid \theta_i,\mathbb{M}_i)p(\theta_i\mid \mathbb{M}_i)d\theta_i,
	\label{eq. marginal likelihood} 
\end{align}
where $\theta_i$ are the parameters of model $\mathbb{M}_i$.

As a point estimate, consider the expectation of $q$:
\vspace{-1.5mm}
\begin{align}
	\mathbb{E}[q\mid \mathcal{D}]
	&=
	\sum_{i=1}^{N} p(\mathbb{M}_i\mid \mathcal{D})\int qp(q\mid \mathbb{M}_i,\mathcal{D})dq.
	&
	\label{eq. point estimate}
\end{align}

The equation \eqref{eq. point estimate} may be viewed as a superposition of individual predictions of the model ensemble, each of which is weighted by $w_i:=p(\mathbb{M}_i\mid \mathcal{D})$, so that
\vspace{-1.5mm}
\begin{align}
	\eqref{eq. point estimate}\ \Leftrightarrow\ 
	\mathbb{E}[q\mid \mathcal{D}]
	&=
	\sum_{i=1}^{N} w_i\mathbb{E}[q\mid \mathbb{M}_i,\mathcal{D}].
	&
	\label{eq. point estimate w_i}
\end{align}

Computing the exact $w_i$ is difficult in general since it involves the marginalization \eqref{eq. marginal likelihood}.
Therefore, BMA may be approximately implemented in practice, e.g., using the Expectation Maximization (EM) algorithm\cite{BMA_forecast_weather, BMA_hydrologic_prediction}, Markov Chain Monte Carlo (MCMC) methods\cite{BMA_long_term_wind_speed}, and Akaike Information Criterion (AIC)-type weighting\cite{BMA_stacking_2018, BMA_mortality_forecasting}.
In this paper, we adopt an AIC-type weighting method, also called pseudo-BMA, which approximates the weight $w_i$ as
\vspace{-1.5mm}
\begin{align}
	w_i=p(\mathbb{M}_i\mid \mathcal{D})
	\approx 
	\cfrac{\exp(\text{elpd}^i)}{\sum_{k=1}^{N}\exp(\text{elpd}^k)},
	\label{eq. approximation of w_i}
\end{align}
where $\text{elpd}^i:=\sum_{j=1}^{n_s}\int p_t(\tilde{q}_j)\log (\tilde{q}_j\mid \mathcal{D},\mathbb{M}_i)d\tilde{q}_j$ denotes the expected log pointwise predictive density of model $\mathbb{M}_i$, where $\{\tilde{q}_j\}_{j=1}^{n_s}$ are unseen new data points and $p_t(\tilde{q}_j)$ is the true distribution of the data. 
In practice, $\text{elpd}^i$ may be also approximated by the Leave-One-Out (LOO) predictor.
For details, refer to \cite{BMA_stacking_2018}.
To compute the weights $w_i$, PyMC\cite{PYMC} is used in the numerical simulations in Section \ref{sec. numerical examples}.

\subsection{Koopman Model Averaging}
Given a data set $\mathcal{D}=\{ (x_i,u_i,y_{i})\mid y_{i}=f(x_i,u_i) \}$, consider an ensemble of $N$ linear embedding models with the common feature maps $z_{k}:=\bm{g}(x_k)\in \mathbb{R}^{N_x}$:
\vspace{-1.5mm}
\begin{numcases}{}
	z^+ = A_i z_k + B_i u_k,
\\
	x^\text{pred}_{k+1} = C_i z^+,
\end{numcases} 
where $i=1,\cdots, N$.
In the proposed algorithm, a base model is trained first by solving Problem \ref{problem 1} on a subset $\mathcal{D}_1$ of the entire date set $\mathcal{D}$ to obtain the feature maps $\bm{g}$ and matrices $(A_1,B_1,C_1)$.
The remaining model parameters $(A_i,B_i,C_i)_{i=2}^N$ are obtained by EDMD, each of which is trained on another subset $\mathcal{D}_i\subset (\mathcal{D}\setminus \mathcal{D}_1)$ of the data:
\vspace{-1.5mm}
\begin{align}
	[A_i\,B_i]&\hspace{-0.5mm}:=\hspace{-0.5mm}
	\underset{[A\ B]}{\text{argmin}}\hspace{-1mm}
	\sum_{(x_j, u_j,y_j)\in \mathcal{D}_{i}}\hspace{-1mm}
	\left\|
		g(y_{j}) - [A\, B]\hspace{-0.5mm}
		\left[
			\begin{array}{c}
				\hspace{-1mm}g(x_j)\hspace{-1mm}
			\\
				u_j
			\end{array}
		\right]
	\right\|_2^2,&
	\label{eq. compute A_i and B_i}
\\
	C_i&:=
	\underset{C}{\text{argmin}}
	\sum_{(x_j, u_j,y_j)\in \mathcal{D}_{i}}
	\left\|
	x_j - 
	Cg(x_{j})
	\right\|_2^2,&
	\label{eq. compute C_i}
\end{align}
where $(A_i,B_i,C_i)$ corresponds to the parameters of model $\mathbb{M}_i$ in Section \ref{sec. BMA outline}.

Given a state $x_k\in \mathcal{X}$ and an input $u_k\in \mathcal{U}$, we assume that both $z_{k+1}=\bm{g}(x_{k+1})$ and $x_{k+1}$ have Gaussian distributions conditioned on the $i$-th model:
\vspace{-2mm}
\begin{align}
	&p(z_{k+1}\mid \mathbb{M}_i,\mathcal{D}) =  
	\mathcal{N}(z_{k+1}; A_i z_k + B_i u_k, \Sigma_z), 
	&
\\
	&
	p(x_{k+1}\mid \mathbb{M}_i,\mathcal{D})\hspace{-0.5mm} =\hspace{-0.5mm}  
	\mathcal{N}(x_{k+1}\hspace{-0.5mm};\hspace{-0.5mm} C_i(A_i z_k + B_i u_k), \Sigma), 
	&
\end{align}
where $\Sigma_z$ and $\Sigma$ are covariance matrices of the distributions.

Taking the expectation \eqref{eq. point estimate w_i} w.r.t. $z_{k+1}$ and $x_{k+1}$, we have:
\begin{align}
	\mathbb{E}[z_{k+1}\mid \mathcal{D}]&=
	\sum_{i=1}^{N} w_i	(A_i z_k + B_i u_k)
	&\nonumber
\\
	&=
	\left(
		\sum_{i=1}^{N} w_iA_i
	\right)z_k 
	+
	\left(
		\sum_{i=1}^{N} w_i B_i
	\right)u_k,
	&
	\label{eq. point estimate of z_{k+1}}
\\
	\mathbb{E}[x_{k+1}\hspace{-0.5mm}\mid\hspace{-0.5mm} \mathcal{D}]&\hspace{-0.5mm}=\hspace{-1mm}
	\left(\hspace{-0.5mm}
	\sum_{i=1}^{N} w_i C_i A_i
	\hspace{-0.5mm}\right)\hspace{-0.9mm}z_k 
	\hspace{-0.5mm}+\hspace{-1mm}
	\left(\hspace{-0.5mm}
	\sum_{i=1}^{N} w_i C_i B_i
	\hspace{-0.5mm}\right)\hspace{-0.9mm}u_k. 
	&
	\label{eq. point estimate of x_{k+1}}
\end{align}

The weight $w_i$ is then approximately computed according to \eqref{eq. approximation of w_i}.
Finally, equations \eqref{eq. point estimate of z_{k+1}} and \eqref{eq. point estimate of x_{k+1}} yield a new \textit{weighted} linear embedding model:
\begin{numcases}{}
	z_{k+1} \approx 
	\left(
	\sum_{i=1}^{N} w_iA_i
	\right)
	z_k
	+
	\left(
	\sum_{i=1}^{N} w_iB_i
	\right)
	u_k, 
	\label{eq. proposed model embedded}
\\
	x^\text{pred}_{k+1}
	=
	\left(
	\sum_{i=1}^{N} w_i C_i A_i
	\right)z_k 
	+
	\left(
	\sum_{i=1}^{N} w_i C_i B_i
	\right)u_k. 
	\label{eq. proposed model state}
\end{numcases}

The outputs of this model are point estimates of the posterior that take the $N$ models' possible outputs into account and are expected to possess a better predictive capability as well as generalizability compared to obtaining a single model only.
The proposed method yields a linear embedding model, which can be used for not only prediction but also other decision-making tasks such as controller designs while it is still model-uncertainty aware by computing the point estimates of the posterior.
The proposed Koopman Model Averaging (KMA) is summarized in Algorithm \ref{algorithm proposed method}.

\begin{algorithm}
	\caption{Koopman Model Averaging (KMA)}
	\renewcommand{\algorithmicrequire}{\textbf{Input:}}
	\renewcommand{\algorithmicensure}{\textbf{Output:}}
	\begin{algorithmic}[1]
		\Require Data set $\mathcal{D}=\{ (x_i,u_i,y_{i})\mid y_{i}=f(x_i,u_i) \}$ 
		\Statex \textbf{Step 1: Training a base model}\\ 
            Solve Problem \ref{problem 1} using a subset $\mathcal{D}_1\subset \mathcal{D}$ of the data to obtain the feature maps $\bm{g}$ and $(A_1,B_1,C_1)$ of the base model
            \vspace{1mm}
            \Statex \textbf{Step 2: Model averaging}
		\For {$i=2:N$}
		\State Compute $(A_i,B_i,C_i)$ in \eqref{eq. compute A_i and B_i} and \eqref{eq. compute C_i} using a subset $\mathcal{D}_i\subset (\mathcal{D}\setminus \mathcal{D}_1)$
		\EndFor 
		\State Extract a subset $\mathcal{D}_a:= (\mathcal{D}\setminus\cup_i \mathcal{D}_i)$ of unseen data for the model and compute $\{w_i\}_{i=1}^N$ in \eqref{eq. approximation of w_i}
		\State Use \eqref{eq. proposed model state} for state prediction
		\State Use \eqref{eq. proposed model embedded} for control application
	\end{algorithmic}
	\label{algorithm proposed method}
\end{algorithm}

\section{Numerical evaluations}
\label{sec. numerical examples}
\subsection{Duffing Oscillator}
\vspace{-1mm}
We first evaluate the proposed model averaging method using the Duffing oscillator with a control input, which is given by:
\vspace{-3mm}
\begin{align}
	\ddot{x}(t)
	=
	-0.5\dot{x}(t) + x(t) - 4 x^3(t) + u(t),
	\label{eq. duffing oscillator}
\end{align}
where the state $x(t)$ and the input $u(t)$ are continuous variables.
It is assumed that the time-series data $x(k\Delta t)$, $u(k\Delta t)$ ($k=0,1,\cdots$) is available to train models, where $\Delta t$ is the sampling period.
Equation \eqref{eq. duffing oscillator} yields a difference equation of the same form as \eqref{eq. target dynamics} with a first-order time discretization so that we can relate the time-series data and the discrete-time dynamics \eqref{eq. target dynamics} as $x(k\Delta t)=x_k$, $u(k\Delta t)=u_k$.
The sampling period is set to $\Delta t=0.01$ in the simulations.

To train a base model, 300 trajectories of the state $x_k$ are generated, each of which has a length of 50 steps and starts from an initial condition sampled from $\text{Uniform}[-3,3]^2$. This corresponds to the data $\mathcal{D}_1$ in Algorithm \ref{algorithm proposed method}. Inputs are sampled from a uniform distribution $u_k\sim \text{Uniform}[-2.5,2.5]$, $\forall k$.
Following a common practice in the literature, we adopt a specific structure of feature maps of the form:
\vspace{-4mm}
\begin{align}
	\bm{g}(x_k)=[x_k^\tr\ g_1(x_k)\ g_2(x_k)\cdots]^\tr,
	\label{eq. feature map with state identity}
\end{align}
where the first $n$ components are the state $x_k$ itself.
This allows an analytic expression of the decoder, i.e., with $C=[I_n\ \bm{0}]$ in \eqref{eq. linear embedding model non autonomous}, we can recover the original state $x_k=C\bm{g}(x_k)$ without learning the parameter $C$.
We have one additional feature map $g_1(x_k)$ in \eqref{eq. feature map with state identity}, which is a neural network with a single hidden layer consisting of 10 neurons.

To populate the model ensemble, we use 4 data subsets $\{ \mathcal{D}_i \}_{i=2}^5$ (so that $N=5$ in Algorithm \ref{algorithm proposed method}). Each $\mathcal{D}_i$ consists of 100 trajectories, each of which is a length of 50 steps.
Data set $\mathcal{D}_a$ consists of 50 trajectories with a length of 20 steps, each of which is sampled from the same distributions of that of $\mathcal{D}_i$.
We use Pytorch to solve Problem \ref{problem 1}. PyMC\cite{PYMC} is used to compute $w_i$ in \eqref{eq. approximation of w_i}.

To compare with the proposed method, we also train two models: an EDMD model and a neural network-based model obtained by Problem \ref{problem 1}.
The neural network-based model is labeled \textit{normal NN model} in the sequel.
Pre-specified feature maps for the EDMD model are monomials up to the second order.
Both the EDMD and the normal NN models are trained on the entire data set $\mathcal{D}=\{\mathcal{D}_i\}_{i=1}^5+\mathcal{D}_a$ used for learning the proposed model.

We consider two control applications in addition to the state prediction: stabilization by LQR and linear MPC. 
The objective of the LQR problem is to have the state $x_k$ converge to $0$. 
In the MPC design, we define a cost function so that the first component of the state $x_k$ will follow a reference signal: 
\vspace{-5mm}
\begin{align}
	r(t)=
	\left\{
		\begin{array}{r}
			-1\ (t\leq 10)
		\\
			1\ (t>10)
		\end{array}
	\right.. 
\end{align}

\begin{figure}[t]
	\centering
	\begin{subfigure}{0.49\linewidth}
		\centering
		\includegraphics[width=\linewidth]{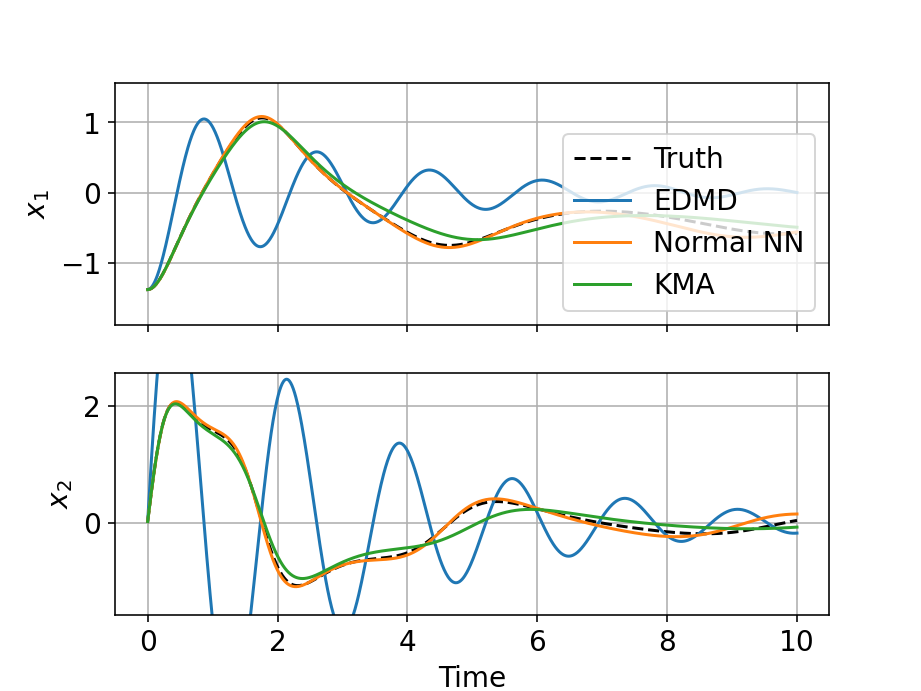}
		\caption{State prediction (IC 1).}
		\label{subfig. duffing state prediction ic 1}
	\end{subfigure}
	\begin{subfigure}{0.49\linewidth}
		\centering
		\includegraphics[width=\linewidth]{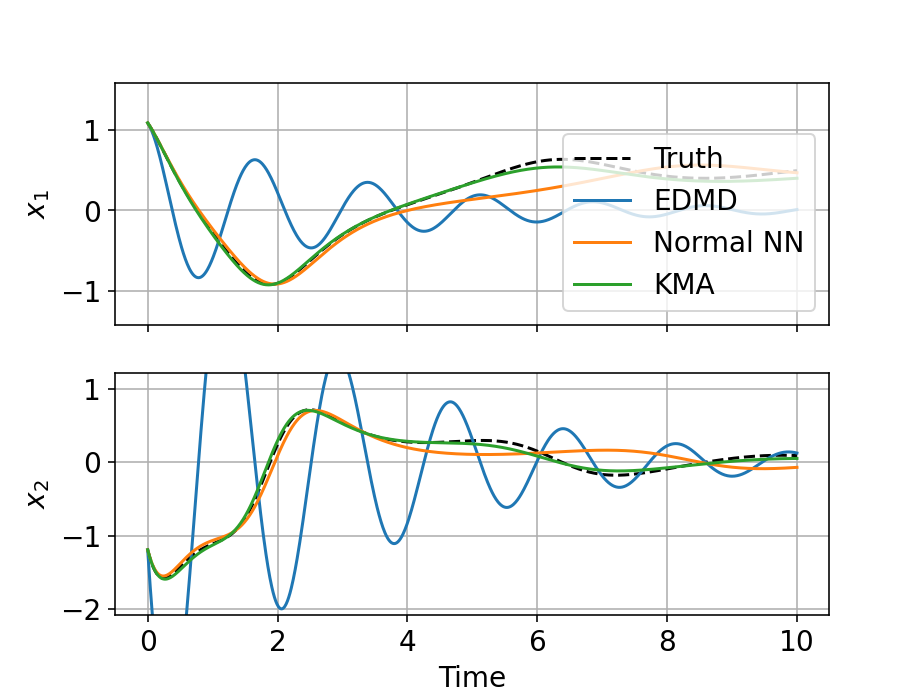}
		\caption{State prediction (IC 2).}
		\label{subfig. duffing state prediction ic 2}
	\end{subfigure}
	\begin{subfigure}{0.49\linewidth}
		\centering
		\includegraphics[width=\linewidth]{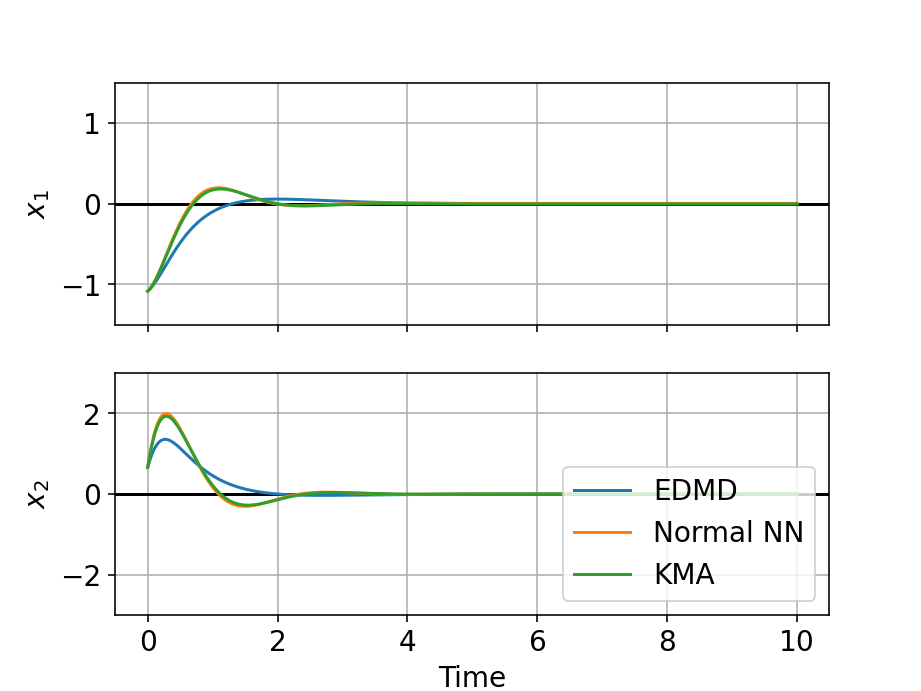}
		\caption{LQR (IC 1).}
		\label{subfig. duffing LQR ic 1}
	\end{subfigure}
	\begin{subfigure}{0.49\linewidth}
		\centering
		\includegraphics[width=\linewidth]{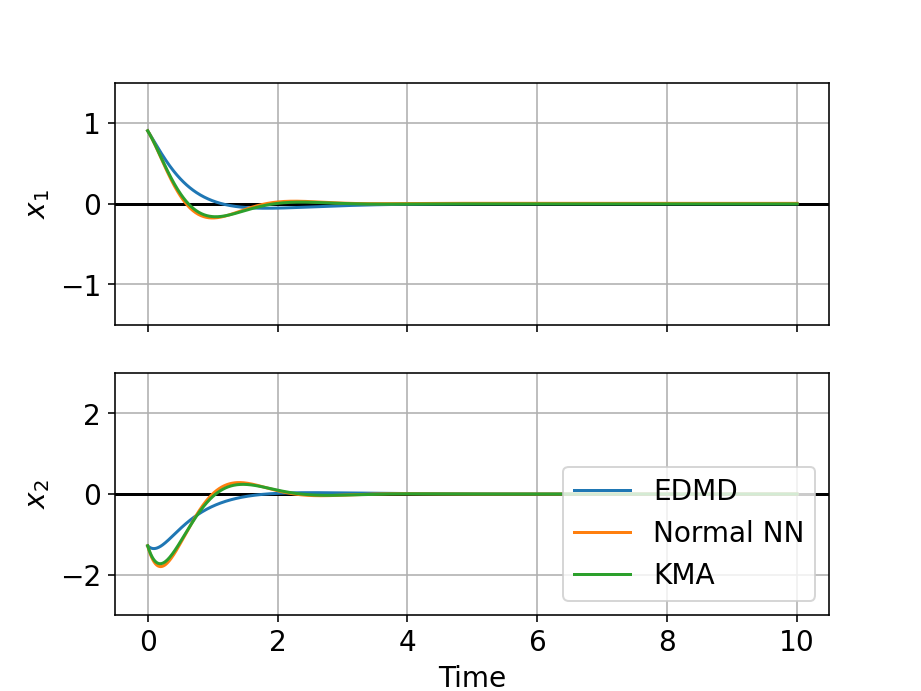}
		\caption{LQR (IC 2).}
		\label{subfig. duffing LQR ic 2}
	\end{subfigure}
	\begin{subfigure}{0.49\linewidth}
		\centering
		\includegraphics[width=\linewidth]{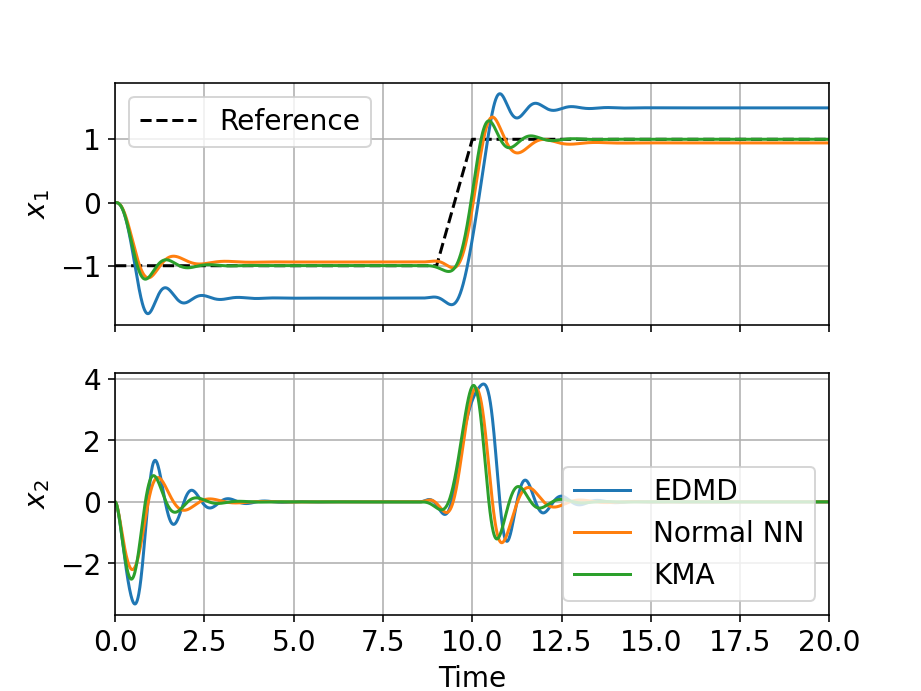}
		\caption{MPC.}
		\label{subfig. duffing MPC}
	\end{subfigure}
	\caption{Duffing oscillator.}
	\label{fig. duffing}
\end{figure}

The results are shown in Fig. \ref{fig. duffing}.
Figures \ref{subfig. duffing state prediction ic 1} and \ref{subfig. duffing state prediction ic 2} are the results of state prediction, where two initial conditions are randomly selected and labeled IC 1 and IC 2, respectively.
Since the EDMD model has a simpler feature map design than the other two neural network-based models, it fails to reconstruct the target dynamics.
On the other hand, both the normal NN and the proposed models show comparable performance and they successfully predict the future states.
Figures \ref{subfig. duffing LQR ic 1} and \ref{subfig. duffing LQR ic 2} are the LQR simulation results with randomly selected initial conditions labeled IC 1 and IC 2, respectively.
All the controllers designed for the three models achieve the control objective in this task.
On the contrary, MPC with the EDMD model fails to track the reference signal as shown in Fig. \ref{subfig. duffing MPC}. 
The normal NN model also has a slight steady state error, whereas the proposed weighted model perfectly tracks the reference signal.
The validation loss of the normal NN model is $7.60\times 10^{-6}$, which is smaller by one order of magnitude than that of the base model of the proposed method, which is $1.81\times 10^{-5}$. 
However, the proposed model outperforms in the MPC task, which implies the effectiveness of aggregating an ensemble of models to find more accurate model outputs.
\vspace{-1.5mm}
 
\begin{figure}[t]
	\centering
	\begin{subfigure}{0.49\linewidth}
		\centering
		\includegraphics[width=\linewidth]{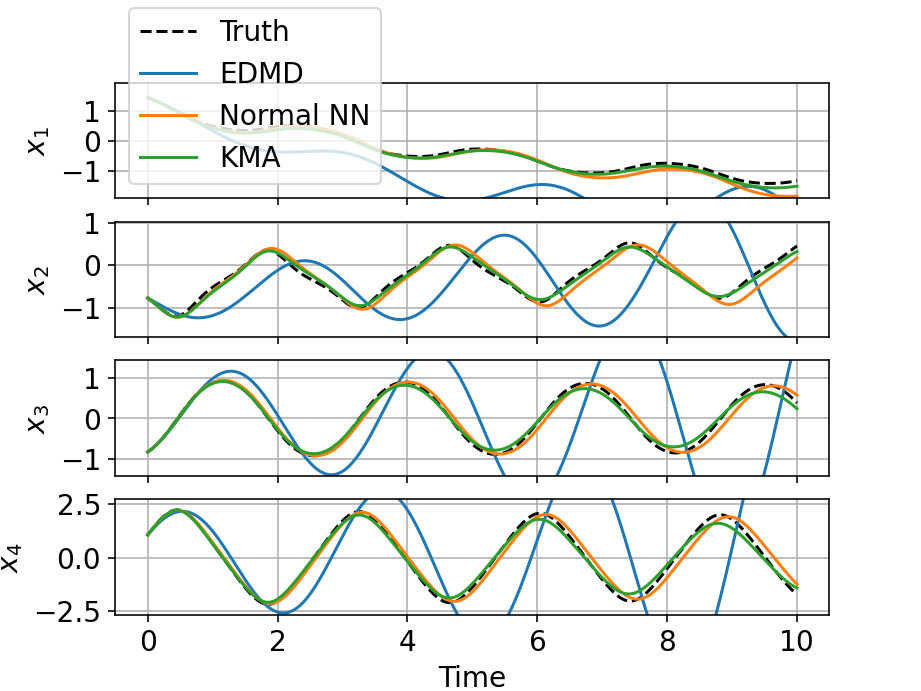}
		\caption{State prediction (IC 1).}
		\label{subfig. cartpole state prediction ic 1}
	\end{subfigure}
	\begin{subfigure}{0.49\linewidth}
		\centering
		\includegraphics[width=\linewidth]{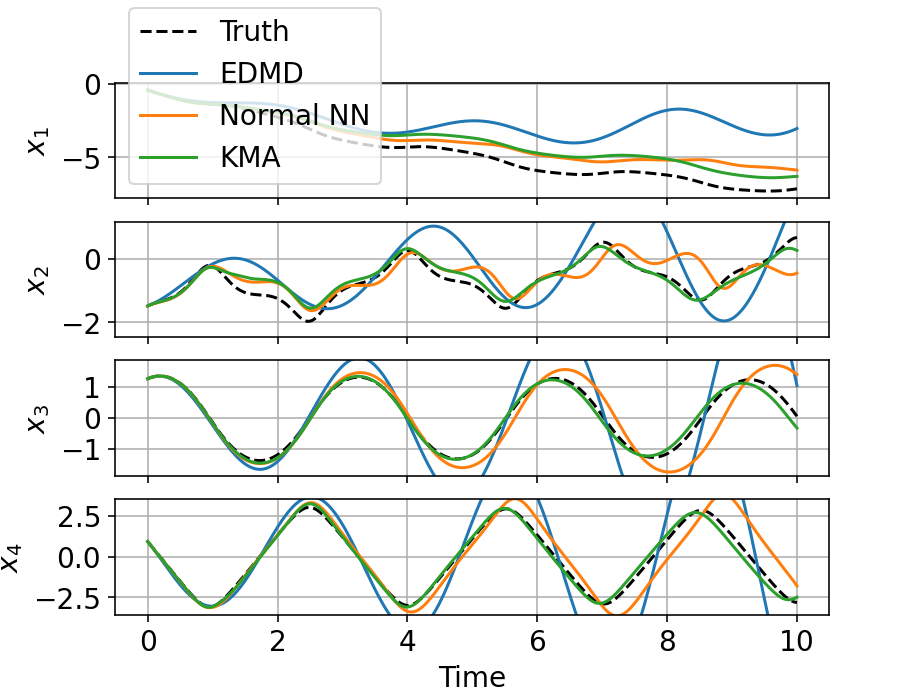}
		\caption{State prediction (IC 2).}
		\label{subfig. cartpole state prediction ic 2}
	\end{subfigure}
	\begin{subfigure}{0.49\linewidth}
		\centering
		\includegraphics[width=\linewidth]{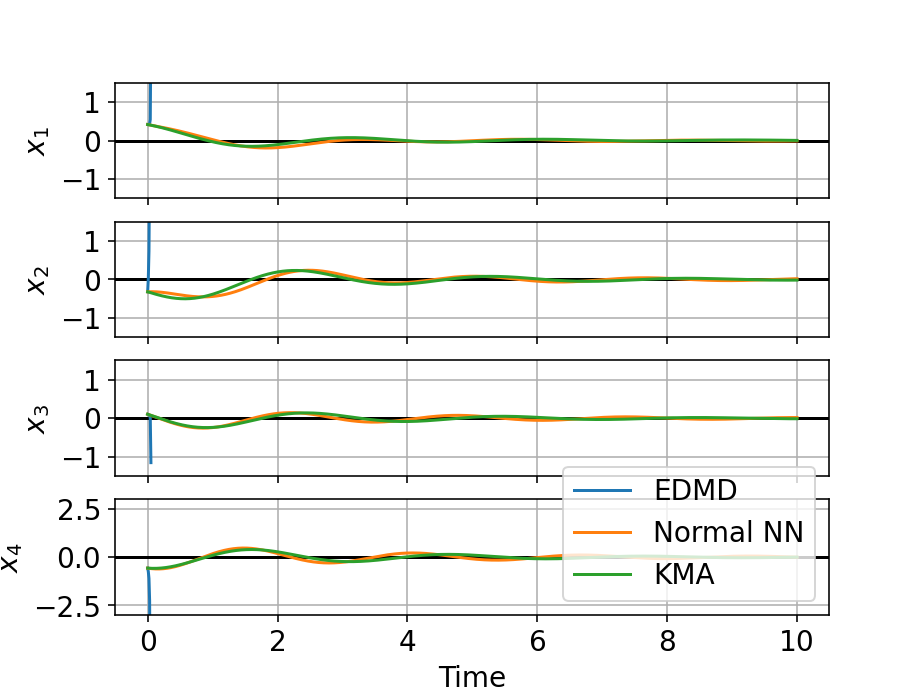}
		\caption{LQR (IC 1).}
		\label{subfig. cartpole LQR ic 1}
	\end{subfigure}
	\begin{subfigure}{0.49\linewidth}
		\centering
		\includegraphics[width=\linewidth]{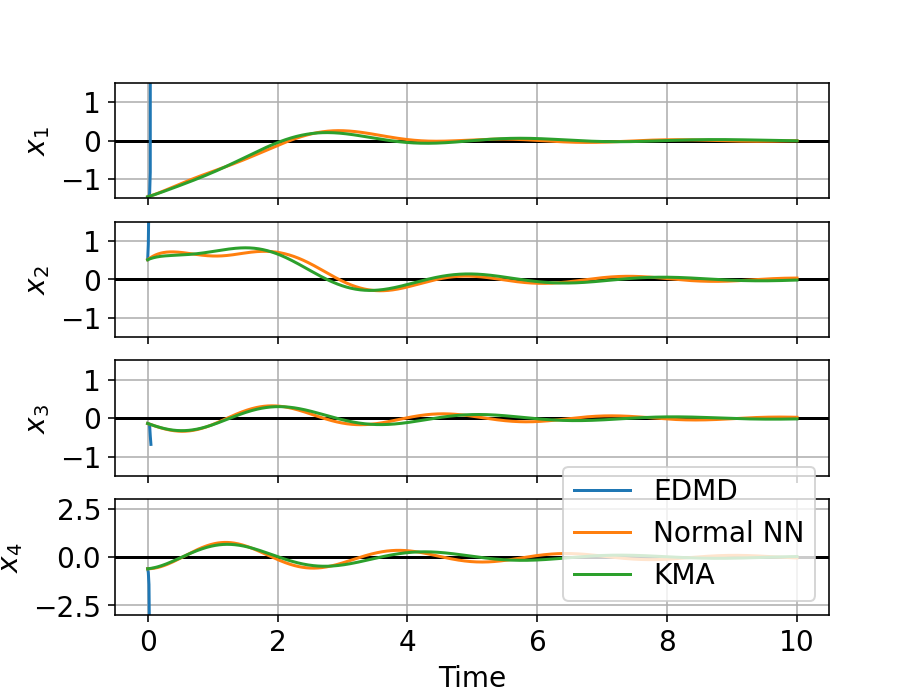}
		\caption{LQR (IC 2).}
		\label{subfig. cartpole LQR ic 2}
	\end{subfigure}
	\begin{subfigure}{0.49\linewidth}
		\centering
		\includegraphics[width=\linewidth]{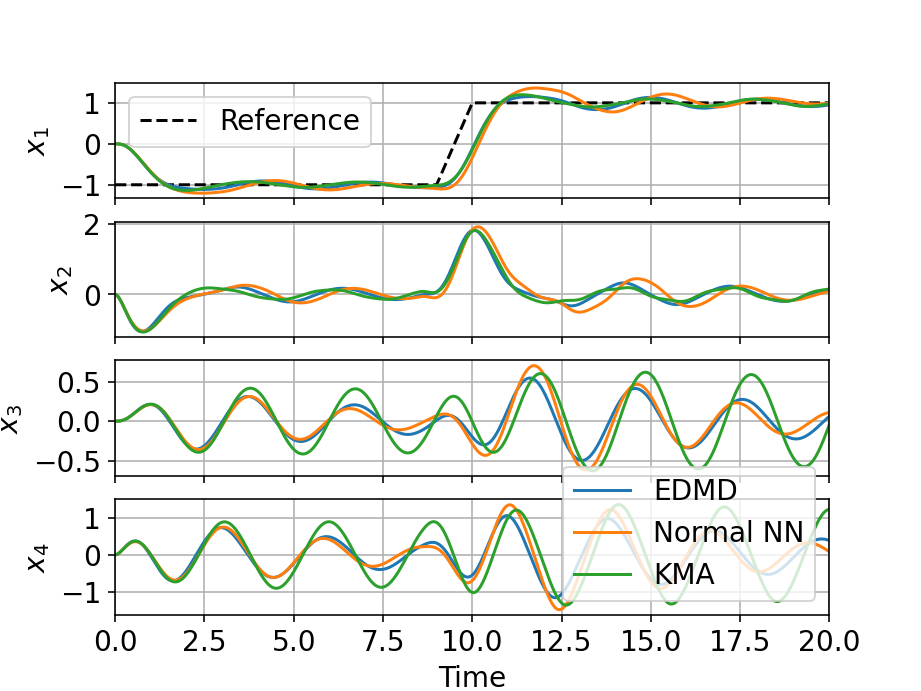}
		\caption{MPC.}
		\label{subfig. cartpole MPC}
	\end{subfigure}
	\caption{Cartpole system.}
	\label{fig. cartpole}
\end{figure}

\subsection{Cartpole}
\vspace{-1mm}
As a more complex nonlinear system, the cartpole is considered as the second example, whose dynamics is given as follows\cite{Data_driven_book}:
\vspace{-2mm}
\begin{align}
	\left\{
	\hspace{-2mm}
	\begin{array}{l}
		\dot{x}_1=x_2,
		\\
		\dot{x}_2=
		\cfrac{
			-\hspace{-1mm}m^2L^2g\cos x_3\sin x_3 
			\hspace{-1mm}+\hspace{-1mm} 
			mL^2A(x_2,x_3,x_4) 
			\hspace{-1mm}+\hspace{-1mm}
			mL^2u
		}
		{mL^2(M+m(1-\cos^2 x_3))},
		\\
		\dot{x}_3=x_4,
		\\
		\dot{x}_4=\hspace{-1mm}
		\cfrac{
			\begin{array}{l}
				(m+M)mgL\sin x_3 -mL\cos x_3 A(x_2,x_3,x_4)
				\\
				\hspace{47mm}+mL\cos x_3 u
			\end{array}
		}
		{mL^2(M+m(1-\cos^2 x_3))},
	\end{array}
	\right.
	\nonumber
\end{align}
where 
$A(x_2,x_3,x_4)=mL{x_4}^2\sin x_3 -\delta x_2$, $m=1$, $M=5$, $L=2$, $g=-10$, and $\delta=1$.
Data collecting procedures and model learning conditions are the same as the first example except for the number of hidden layers of the neural network, which is two for the cartpole system.

The results are shown in Fig. \ref{fig. cartpole}. 
In this example, all the models are successfully applied to the MPC task as in Fig. \ref{subfig. cartpole MPC}. 
On the other hand, the EDMD model fails in the LQR design problem (Figs. \ref{subfig. cartpole LQR ic 1} and \ref{subfig. cartpole LQR ic 2}), which is considered as a result of too simple feature map design for the four dimensional dynamics of the cartpole.
The results of the state prediction are shown in Figs. \ref{subfig. cartpole state prediction ic 1} and \ref{subfig. cartpole state prediction ic 2}. The EDMD model also has difficulty in this task. 
Both the normal NN and the proposed models show reasonable predictions with the initial condition IC 1 (Fig. \ref{subfig. cartpole state prediction ic 1}). 
However, the prediction of the normal NN model starts deviating from the true values at around $t=5$ with the second initial condition (Fig. \ref{subfig. cartpole state prediction ic 2}). 
The proposed model, on the other hand, shows better accuracy and its prediction follows the true values until the end of the simulation.
It is noted that the validation loss of the normal NN model ($1.16\times 10^{-5}$) is smaller than that of the proposed method ($2.82\times 10^{-5}$). This result implies that the proposed model effectively combines the model ensemble into a new weighted model to acquire high predictive accuracy with respect to wider range of regimes of dynamics than the normal NN model.
\vspace{-2mm}

\section{CONCLUSION}
\vspace{-2mm}
We propose a model averaging method for learning Koopman operator-based models for prediction and control utilizing the Bayesian model averaging. 
While the Koopman operator framework allows to obtain linear embedding models from data so that linear systems theories can be applied to control possibly nonlinear dynamics, it is challenging to accurately reconstruct the dynamics and deploy the learned model in several applications. 
Considering that training a single model only may not be sufficient to obtain model predictions that are accurate enough for a wide range of operating points even if the model possesses high accuracy w.r.t. a certain regime of dynamics such as training data, the proposed method first trains a base model with the use of neural networks and then populates an ensemble of models on different data points. 
Based on the ideas of the Bayesian model averaging, these models are merged into a new weighted linear embedding model, in which individual outputs of the model ensemble are weighted according to the posterior model evidence.
This model explicitly takes  uncertainty into account and the overall performance is expected to be improved.
Using numerical simulations, it is shown that the proposed weighted model achieves better state predictive accuracy as well as greater generalizability to different control applications compared to other Koopman operator-based models.
\vspace{-5mm}
\bibliography{ACC2025}  
\bibliographystyle{ieeetr}  

\end{document}